\documentclass[11pt]{amsart}
\usepackage{amsmath,amssymb,mathtools}
\usepackage[hidelinks]{hyperref}
\newtheorem{theorem}{Theorem}
\newtheorem{lemma}[theorem]{Lemma}

\theoremstyle{remark}
\newtheorem{remark}[theorem]{Remark}
\newcommand{\T}{\mathbb T}
\newcommand{\D}{\mathbb D}
\newcommand{\C}{\mathbb C}
\newcommand{\Z}{\mathbb Z}
\newcommand{\N}{\mathbb N}
\newcommand{\Pp}{P_{+}}
\newcommand{\norm}[2]{\lVert #1\rVert_{L_{#2}}}
\newcommand{\dm}{\,dm}
\newcommand{\dA}{\,dA}
\numberwithin{equation}{section}
\title[Extremizers for Riesz projections]{Extremizers for contractive Riesz projections}
\date{November 26, 2026}
\author[L.~M.~Casta\~no-Mar\'in]{Luis Miguel Casta\~no-Mar\'in}
\address{Departamento de Matem\'aticas\\
Universidad Nacional de Colombia\\
Bogot\'a, Colombia}
\email{lucastanom@unal.edu.co}
\subjclass[2020]{42B05, 30H10, 42B30}
\keywords{Riesz projection, Hardy spaces, inner functions, extremizers, subharmonic functions}
\begin{document}
\begin{abstract}
Brevig, Llinares and Seip conjectured that the nonzero extremizers of the contractive inequality $\|P_+f\|_{L_p}\leq\|f\|_{L_q}$, where $p=4(1-1/q)$ and $1<q\leq\infty$, $q\ne2$, are precisely the nonzero constant multiples of inner functions. We prove this conjecture by applying Green's formula to the subharmonic minorant of Fu and Li. The one-variable characterization also yields the equality cases in the corresponding iterated contractive inequalities on the polytorus.
\end{abstract}
\maketitle

\section{Introduction}
Throughout the paper, $\D$ denotes the open unit disc, $\T=\partial\D$ and $m$ is normalized Lebesgue measure on $\T$. For a positive integer $d$, we write $m_d=m^{\otimes d}$ and, for a measurable function $f:\T^d\to\C$, set
\begin{equation}\label{eq:norm}
 \norm{f}{s}=\left(\int_{\T^d}|f|^s\,dm_d\right)^{1/s}
 \quad(0<s<\infty),\qquad
 \norm{f}{\infty}=\operatorname*{ess\,sup}_{z\in\T^d}|f(z)|.
\end{equation}
As usual, $L_s(\T^d)$ consists of the functions for which this quantity is finite, identified up to equality almost everywhere; for $0<s<1$ it is a quasi-norm.

For $f\in L_1(\T^d)$ and $\alpha=(\alpha_1,\ldots,\alpha_d)\in\Z^d$, let
\begin{equation}\label{eq:fourier}
 \widehat f(\alpha)=\int_{\T^d} f(z)\overline{z^\alpha}\,dm_d(z),
 \qquad z^\alpha=\prod_{j=1}^d z_j^{\alpha_j}.
\end{equation}
Set $\N_0=\{0,1,2,\ldots\}$. For $1\leq s\leq\infty$, the Hardy space $H_s(\T^d)$ is the subspace of $L_s(\T^d)$ given by
\begin{equation}\label{eq:hardy}
 H_s(\T^d)=\{f\in L_s(\T^d):\widehat f(\alpha)=0
 \text{ for }\alpha\notin\N_0^d\}.
\end{equation}
We identify its elements with their holomorphic extensions to $\D^d$. An inner function on $\D^d$ is a bounded holomorphic function whose radial boundary values have modulus one almost everywhere on $\T^d$.

The Riesz projection is the orthogonal projection $\Pp^{(d)}:L_2(\T^d)\to H_2(\T^d)$. On trigonometric polynomials its action is given by
\begin{equation}\label{eq:projection}
 \Pp^{(d)}f(z)=\sum_{\alpha\in\N_0^d}\widehat f(\alpha)z^\alpha.
\end{equation}
For $1<q<\infty$, it extends to a bounded projection from $L_q(\T^d)$ onto $H_q(\T^d)$, by the classical theorem of M. Riesz and Fubini's theorem. On $L_\infty(\T^d)$ we use the restriction of the $L_2$ projection. We write $\Pp=\Pp^{(1)}$.

Marzo and Seip \cite{MS} initiated the study of the critical exponent of the Riesz projection. In dimension one and for $1<q\leq\infty$, this exponent is
\begin{equation}\label{eq:critical}
 \mathfrak p_1(q)=\sup\big\{s>0:\norm{\Pp f}{s}\leq\norm{f}{q}
 \text{ for every } f\in L_q(\T)\big\}.
\end{equation}
They proved that $\mathfrak p_1(\infty)=4$. Brevig, Ortega-Cerd\`a, Seip and Zhao \cite[Conjecture 4]{BOSZ} conjectured that
\begin{equation}\label{eq:criticalvalue}
 \mathfrak p_1(q)=4\left(1-\frac1q\right)\quad(1<q\leq\infty).
\end{equation}
Fu and Li \cite[Theorem 1.1]{FL} recently proved \eqref{eq:criticalvalue} for all finite $q>1$.

In \cite[Theorem 6 and Corollary 7]{BLS}, Brevig, Llinares and Seip proved that the nonzero extremizers at $q=\infty$ and $q=4/3$ are exactly the nonzero scalar multiples of inner functions. They conjectured that the same characterization holds for every $1<q\leq\infty$, $q\ne2$ \cite[Section 3]{BLS}. At $q=2$, orthogonality gives equality for every nonzero element of $H_2(\T)$. The remaining cases of their conjecture are settled by the following theorem.

\begin{theorem}\label{thm:one}
Let $1<q\leq\infty$, $q\ne2$, and $p=4(1-1/q)$. If $f\in L_q(\T)$ is nonzero, then
\begin{equation}\label{eq:equalityone}
 \norm{\Pp f}{p}=\norm{f}{q}
\end{equation}
if and only if $f$ is a constant multiple of an inner function.
\end{theorem}

The subharmonic minorant constructed by Fu and Li \cite{FL} has a Laplacian that is strictly positive away from the origin. Green's formula turns this property into an integral lower bound for the deficit in the contractive inequality. Analyticity and strict concavity then yield the equality cases. The same argument applies when $0<p<1$.

For $d\geq1$, define
\begin{equation}\label{eq:ad}
 a_d(q)=2+\frac{2(q-2)}{2+d(q-2)}
 \quad\left(\frac{2d}{d+1}<q<\infty\right),\qquad
 a_d(\infty)=2+\frac{2}{d}.
\end{equation}
Brevig, Llinares and Seip \cite[Section 1]{BLS} conjectured that $a_d(q)$ is the critical exponent in this range. Their iteration inequality \cite[Lemma 10]{BLS}, together with the one-variable theorem of Fu and Li, yields contractivity at $a_d(q)$; for $d=2$, optimality follows from \cite[Theorem 3]{BLS}. Applying Theorem~\ref{thm:one} to the coordinate projections yields the equality cases for every $d$.

\begin{theorem}\label{thm:many}
Let $d\geq1$ and $2d/(d+1)<q\leq\infty$, $q\ne2$. If $f\in L_q(\T^d)$ is nonzero, then
\begin{equation}\label{eq:equalitymany}
 \norm{\Pp^{(d)}f}{a_d(q)}=\norm{f}{q}
\end{equation}
if and only if $f$ is a constant multiple of an inner function on $\D^d$.
\end{theorem}

\section{Proof of Theorem \ref{thm:one}}
Fix $1<q<\infty$, $q\ne2$, and put
\begin{equation}\label{eq:parameters}
 p=4(1-1/q),\qquad \lambda=\frac pq,\qquad P=\frac p2,\qquad Q=\frac q2.
\end{equation}
In particular,
\begin{equation}\label{eq:relations}
 q-p=\frac{(q-2)^2}{q}>0,\qquad
 PQ=2Q-1>0,\qquad Q-P=\frac{(Q-1)^2}{Q}>0.
\end{equation}
Thus $0<\lambda<1$. Fu and Li \cite[Lemma 2.1]{FL} introduced the subharmonic minorant $G:\C\to\mathbb R$ given by
\begin{equation}\label{eq:G}
 G(\zeta)=\inf_{t>0}F(t,\zeta),\qquad
 F(t,\zeta)=\frac{P}{Q}\frac{|t+\zeta|^{2Q}}{t^Q}+1-\frac{P}{Q}-t^P.
\end{equation}
The infimum is finite by \cite[Lemma 3.1]{FL}, and the strict positivity of its Laplacian away from zero follows from the calculation in \cite[Lemma 2.1]{FL}.

\begin{lemma}[Fu and Li \cite{FL}]\label{lem:G}
The function $G$ in \eqref{eq:G} is of class $C^\infty$, $G(0)=0$, and its Euclidean Laplacian satisfies $\Delta G\geq0$ on $\C$ and $\Delta G(\zeta)>0$ for every $\zeta\ne0$. Moreover, for $u,v\in\C$,
\begin{equation}\label{eq:minorant}
 \lambda|u+\overline v|^q+1-\lambda-|u|^p\geq G(uv).
\end{equation}
\end{lemma}
\begin{proof}
By \cite[Lemma 3.1]{FL}, the infimum in \eqref{eq:G} is attained at a unique $t=t(\zeta)>|\zeta|$, depending smoothly on $\zeta$, and
\begin{equation}\label{eq:t}
 F_{tt}(t(\zeta),\zeta)>0.
\end{equation}
We have $G(0)=0$ since the weighted arithmetic--geometric mean inequality gives
\[
 F(t,0)=\lambda t^Q+1-\lambda-t^P\geq0=F(1,0).
\]
The calculation in the proof of \cite[Lemma 2.1]{FL} gives
\begin{equation}\label{eq:lap}
 \Delta G(\zeta)=\frac{4P^3Q\,t(\zeta)^{2P-2}|\zeta|^2}
 {\big(t(\zeta)^2-|\zeta|^2\big)^2F_{tt}(t(\zeta),\zeta)}.
\end{equation}
The denominator is strictly positive, so $\Delta G(\zeta)>0$ for $\zeta\ne0$ and $\Delta G(0)=0$.

If $u\ne0$, then
\[
 |u+\overline v|^q=\frac{\big||u|^2+uv\big|^q}{|u|^q}.
\]
Consequently, the left-hand side of \eqref{eq:minorant} is $F(|u|^2,uv)$, which is bounded below by $G(uv)$. If $u=0$, it is $\lambda|v|^q+1-\lambda\geq0=G(0)$.
\end{proof}

Green's formula converts the strict positivity in \eqref{eq:lap} into an estimate for the deficit. Let $dA$ denote unnormalized planar Lebesgue measure.
\begin{lemma}\label{lem:deficit}
Let $f\in L_q(\T)$, $\norm{f}{q}=1$, and write $f=g+\overline h$, where $g=\Pp f\in H_q(\T)$ and $h\in H_q(\T)$ satisfies $h(0)=0$. With $G$ as in \eqref{eq:G},
\begin{equation}\label{eq:deficit}
 1-\norm{g}{p}^p\geq\frac{1}{2\pi}\int_{\D}
 \log\frac{1}{|z|}\,\Delta G(g(z)h(z))\,|(gh)'(z)|^2\dA(z).
\end{equation}
\end{lemma}
\begin{proof}
The classical Hardy space decomposition $f=g+\overline h$ follows from the boundedness of $\Pp$ on $L_q(\T)$, with $h(0)=0$ since $f-g$ has only strictly negative Fourier coefficients. For $0<r<1$ set $g_r(z)=g(rz)$, $h_r(z)=h(rz)$ and $f_r=g_r+\overline{h_r}$ on $\T$. The function $U=G\circ(gh)$ is smooth on $\D$, and $U(0)=G(0)=0$. Integrating \eqref{eq:minorant} on $\T$ gives
\begin{equation}\label{eq:radius}
 \lambda\norm{f_r}{q}^q+1-\lambda-\norm{g_r}{p}^p
 \geq\int_{\T}U(rz)\dm(z).
\end{equation}
Since $gh$ is holomorphic, the chain rule for the Laplacian gives
\begin{equation}\label{eq:chain}
 \Delta U(z)=\Delta G(g(z)h(z))\,|(gh)'(z)|^2\geq0.
\end{equation}
The classical Green formula on the disc of radius $r$, with $U(0)=0$, gives
\begin{equation}\label{eq:green}
 \int_{\T}U(rz)\dm(z)=\frac{1}{2\pi}\int_{|z|<r}
 \log\frac r{|z|}\,\Delta U(z)\dA(z).
\end{equation}
By the classical radial convergence theorem in $H_q$, $g_r\to g$ and $h_r\to h$ in $L_q(\T)$, hence $f_r\to f$ in $L_q(\T)$. Since $p<q$, we also have $g_r\to g$ in $L_p(\T)$. When $p<1$, subadditivity gives
\[
 \big|\norm{g_r}{p}^p-\norm{g}{p}^p\big|
 \leq\int_{\T}|g_r-g|^p\dm\leq\norm{g_r-g}{q}^p.
\]
It follows that the left-hand side of \eqref{eq:radius} converges to $1-\norm{g}{p}^p$. On the right-hand side of \eqref{eq:green}, the nonnegative integrands, extended by zero outside $\{|z|<r\}$, increase pointwise almost everywhere as $r$ increases to one. Monotone convergence, \eqref{eq:chain}, and \eqref{eq:radius} give \eqref{eq:deficit}.
\end{proof}

\begin{proof}[Proof of Theorem \ref{thm:one}]
The case $q=\infty$ is \cite[Theorem 6]{BLS}. Suppose $1<q<\infty$ and normalize $f$ so that $\norm{f}{q}=1$. If $\norm{\Pp f}{p}=1$, the left-hand side of \eqref{eq:deficit} is zero. By Lemma \ref{lem:G}, the integrand is nonnegative. Since $\log(1/|z|)>0$ for $0<|z|<1$ and the function $z\mapsto\Delta G(g(z)h(z))|(gh)'(z)|^2$ is continuous, \eqref{eq:deficit} implies
\begin{equation}\label{eq:vanish}
 \Delta G(gh)\,|(gh)'|^2=0\qquad\text{on }\D.
\end{equation}
On the open set where $gh\ne0$, Lemma \ref{lem:G} and \eqref{eq:vanish} give $(gh)'=0$. If that set is nonempty, the identity theorem makes $gh$ constant on $\D$. Since $h(0)=0$, it follows that $gh\equiv0$. As $\norm{g}{p}=1$, the identity theorem gives $h\equiv0$.

Thus $f=g$ and $\norm{f}{p}=\norm{f}{q}=1$. Since $0<p/q<1$, strict concavity of $s\mapsto s^{p/q}$ on $[0,\infty)$ gives
\[
 1=\int_{\T}(|f|^q)^{p/q}\dm
 \leq\left(\int_{\T}|f|^q\dm\right)^{p/q}=1,
\]
with equality only if $|f|^q$ is constant almost everywhere. The normalization gives $|f|=1$ almost everywhere. The classical Poisson representation of $H_q(\T)$ functions now gives $|f(z)|\leq1$ on $\D$, so $f$ is inner. Removing the normalization yields the claimed form. Conversely, if $f=cI$ with $I$ inner, then $\Pp f=f$ and $\norm{f}{p}=|c|=\norm{f}{q}$.
\end{proof}

\section{Proof of Theorem \ref{thm:many}}
The contractive inequality on the polytorus is obtained by the iteration of Brevig, Llinares and Seip \cite[Lemma 10]{BLS}. The equality cases follow by applying Theorem~\ref{thm:one} to the first projection in each possible order.

For $1< s\leq\infty$, let $b(s)=4(1-1/s)$, with $1/\infty=0$. Fix $q$ and $d$ as in Theorem \ref{thm:many}, and put $r_0=q$, $r_k=a_k(q)$ for $1\leq k\leq d$. The exponents satisfy
\begin{equation}\label{eq:iterate}
 r_{k+1}=b(r_k),\qquad 0\leq k<d.
\end{equation}
Indeed, for finite $q$ and $v=q-2$, $r_k=2+2v/(2+kv)$, so
\[
 b(r_k)=4-\frac{2(2+kv)}{2+(k+1)v}
 =2+\frac{2v}{2+(k+1)v}.
\]
For $q=\infty$, the first iterate is $4$, and for $k\geq1$ the same identity follows from $r_k=2+2/k$. For finite $q<2$, the assumption $q>2d/(d+1)$ gives $2+(d+1)v>0$. Thus $2+kv>0$ for $0\leq k\leq d$, and $r_k>1$ for $0\leq k<d$, since
\[
 r_k-1=\frac{2+(k+2)v}{2+kv}>0\qquad (0\leq k<d).
\]
Also $r_d>0$, since $r_d=2(2+(d+1)v)/(2+dv)$. For $q>2$ all iterates exceed two. None of the iterates is two, and $b(s)<s$ for finite $s\ne2$ by \eqref{eq:relations}. In particular, $r_{k+1}<r_k$ in every case.

For $1\leq j\leq d$, let $R_j:L_s(\T^d)\to L_s(\T^d)$ ($1<s<\infty$) denote the coordinate Riesz projection, specified by
\begin{equation}\label{eq:coord}
 \widehat{R_j F}(\alpha)=\begin{cases}
 \widehat F(\alpha),&\alpha_j\geq0,\\
 0,&\alpha_j<0.
 \end{cases}
\end{equation}
For $s=\infty$ we restrict $R_j$ on $L_2(\T^d)$ to $L_\infty(\T^d)$. On trigonometric polynomials, the $R_j$ commute and $\Pp^{(d)}=R_d\cdots R_1$. By boundedness on $L_s(\T^d)$ for $1<s<\infty$, which follows from the one-dimensional boundedness and Fubini's theorem, these identities extend to $L_s(\T^d)$. For $f\in L_\infty(\T^d)$ the same identities are understood on $L_2(\T^d)$.

\begin{lemma}\label{lem:slice}
Suppose $d\geq2$, $1<s\leq\infty$, $s\ne2$, and $F\in L_s(\T^d)\setminus\{0\}$. Then
\begin{equation}\label{eq:slice}
 \norm{R_j F}{b(s)}\leq\norm{F}{s}.
\end{equation}
If equality holds in \eqref{eq:slice}, then $R_jF=F$ and $|F|=\norm{F}{s}$ almost everywhere.
\end{lemma}
\begin{proof}
Write $(z,y)\in\T\times\T^{d-1}$ for the variables, where $z$ denotes the $j$th coordinate. For almost every $y$, let $F_y(z)=F(z,y)$ and set $A(y)=\norm{F_y}{s}$. The function $A$ is measurable by Fubini's theorem (applied to level sets when $s=\infty$). For almost every $y$ we have
\begin{equation}\label{eq:sections}
 (R_j F)(\cdot,y)=\Pp(F_y).
\end{equation}
This is the usual slice representation of the coordinate projection: it follows for trigonometric polynomials directly and then in $L_s$ by density and boundedness of the one-variable projection. For $s=\infty$ the same argument applies in $L_2$. Put $t=b(s)$. The one-dimensional contractive inequality and Tonelli's theorem give
\begin{equation}\label{eq:slicechain}
 \begin{split}
 \norm{R_jF}{t}^t
 &=\int_{\T^{d-1}}\norm{\Pp F_y}{t}^t\,dm_{d-1}(y)\\
 &\leq\int_{\T^{d-1}}A(y)^t\,dm_{d-1}(y)\leq\norm{F}{s}^t.
 \end{split}
\end{equation}
The last inequality is Jensen's inequality for the concave function $u\mapsto u^{t/s}$ when $s<\infty$, and follows from $A(y)\leq\norm{F}{\infty}$ almost everywhere when $s=\infty$.

Suppose now that equality holds in \eqref{eq:slice}. Equality in the last inequality of \eqref{eq:slicechain} forces $A(y)=\norm{F}{s}>0$ almost everywhere, by strict concavity when $s<\infty$ and by integration of a nonnegative difference when $s=\infty$. Equality in the first inequality of \eqref{eq:slicechain} forces
\[
 \norm{\Pp F_y}{t}^t=A(y)^t
\]
for almost every $y$, again by integration of a nonnegative difference. Theorem \ref{thm:one} applies to each such slice, so $\Pp F_y=F_y$ and $|F_y|=A(y)$ almost everywhere on $\T$. Fubini's theorem gives $R_jF=F$ and $|F|=\norm{F}{s}$ almost everywhere on $\T^d$.
\end{proof}

\begin{proof}[Proof of Theorem \ref{thm:many}]
The case $d=1$ is Theorem \ref{thm:one}. Suppose $d\geq2$. Fix a permutation $\sigma$ of $\{1,\ldots,d\}$ and put
\[
 F_0=f,\qquad F_k=R_{\sigma(k)}F_{k-1}\quad(1\leq k\leq d).
\]
All these projections may be taken in $L_q$ when $q<\infty$ and in $L_2$ when $q=\infty$. If $q=\infty$, they also belong to $L_s$ for every finite $s>1$ by the boundedness on these spaces and uniqueness of Fourier coefficients. In particular, the identities defining $F_k$ are compatible with the slice projections at each intermediate exponent. By Lemma \ref{lem:slice}, \eqref{eq:iterate}, and the commutation of the coordinate projections,
\begin{equation}\label{eq:globalchain}
 \norm{\Pp^{(d)}f}{r_d}=\norm{F_d}{r_d}
 \leq\norm{F_{d-1}}{r_{d-1}}\leq\cdots\leq\norm{F_1}{r_1}\leq\norm{f}{q}.
\end{equation}
The inequality for a zero intermediate function follows from linearity.

Assume equality between the first and last terms of \eqref{eq:globalchain}. Since $\norm{f}{q}>0$, all the intermediate inequalities are equalities. In particular,
\[
 \norm{R_{\sigma(1)} f}{r_1}=\norm{f}{q}.
\]
Lemma \ref{lem:slice} gives $R_{\sigma(1)}f=f$ and $|f|=\norm{f}{q}$ almost everywhere. For each $j$ choose a permutation with $\sigma(1)=j$. We conclude that $R_j f=f$ for every $j$, and consequently $\widehat f(\alpha)=0$ whenever at least one coordinate of $\alpha$ is negative. Also $f\in L_\infty(\T^d)$ since $|f|=\norm{f}{q}$ almost everywhere.

The function $I=f/\norm{f}{q}$ therefore belongs to $H_\infty(\T^d)$ and has modulus one almost everywhere. By the classical identification of $H_\infty(\T^d)$ with bounded holomorphic functions on the polydisc through the product Poisson kernel, $I$ is inner. Conversely, if $f=cI$ with $I$ inner, then $\Pp^{(d)}f=f$ and its $L_s$ norm is $|c|$ for every $0<s\leq\infty$. In particular, equality holds in \eqref{eq:globalchain}.
\end{proof}

\begin{remark}
At $q=2$, orthogonality gives equality precisely for the nonzero elements of $H_2(\T^d)$.
\end{remark}

\end{document}